\documentclass[a4paper, 10 pt, reqno]{article}
\usepackage{amsfonts, amssymb, amsmath, eucal, amsthm}
\usepackage[all]{xy}

\newtheorem{lemma}{Lemma}
\newtheorem{theorem}[lemma]{Theorem}
\newtheorem*{theorem*}{Theorem}
\newtheorem{proposition}[lemma]{Proposition}
\newtheorem{corollary}[lemma]{Corollary}

\theoremstyle{definition}
\newtheorem{definition}[lemma]{Definition}
\newtheorem{note}[lemma]{Note}

\newtheorem{remark}[lemma]{Remark}

\newcommand{\Q}{\mathbb{Q}}
\renewcommand{\H}{\mathbb{H}}
\newcommand{\ev}{\mathrm{ev}}
\newcommand{\Z}{\mathbb{Z}}
\newcommand{\Id}{\mathrm{Id}}
\newcommand{\Th}{\mathrm{Th}}

\newcommand{\ChasSullivan}{ChasSullivan}
\newcommand{\CohenJones}{MR1942249}
\newcommand{\TamanoiBV}{MR2211159}
\newcommand{\Menichi}{MR2466078}

\newcommand{\SpectralSequences}{MR1793722}
\newcommand{\Elements}{MR516508}
\newcommand{\Bott}{MR0102803}
\newcommand{\EvenBott}{MR0087035}
\newcommand{\Hatcher}{MR1867354}
\newcommand{\Godin}{Godin}
\newcommand{\Yang}{Yang}
\newcommand{\Vaintrob}{Vaintrob}

\title{String Topology for Lie Groups}
\author{Richard Hepworth\footnote{The author is supported by E.P.S.R.C.~Postdoctoral Research Fellowship EP/D066980.}\\  Department of Pure Mathematics\\ University of Sheffield}
\date{}

\begin{document}
\maketitle

\begin{abstract}
\noindent In 1999 Chas and Sullivan showed that the homology of the free loop space of an oriented manifold admits the structure of a Batalin-Vilkovisky algebra.  In this paper we give a direct description of this Batalin-Vilkovisky algebra in the case that the manifold is a compact Lie group $G$.  Our answer is phrased in terms of the homology of $G$, the homology of the space of based loops on $G$, and the homology suspension.  The result is applied to compute the Batalin-Vilkovisky algebra associated to the special orthogonal groups $SO(n)$ with coefficients in the rational numbers and in the integers modulo two.
\end{abstract}

\section{Introduction}

In their seminal paper \cite{\ChasSullivan}, Chas and Sullivan introduced two new algebraic operations on the homology groups $\H_\ast(LM)=H_{\ast+\dim M}(LM)$ of the free loop space of a closed oriented manifold $M$.  The first of these operations is the \emph{loop product}
\[\H_\ast(LM)\otimes\H_\ast(LM)\to\H_\ast(LM),\]
whose definition combines the intersection of cycles in $M$ and the concatenation of loops in $M$.  It makes $\H_\ast(LM)$ into an associative, commutative graded ring with unit.  The second of these operations is the \emph{BV-operator}
\[\Delta\colon\H_\ast(LM)\to\H_{\ast+1}(LM)\]
obtained by rotating loops.  It satisfies $\Delta^2=0$ and also the identity
\begin{eqnarray*}\Delta(abc)
&=&\Delta(ab)c+(-1)^{|a|}a\Delta(bc)+(-1)^{(|a|-1)|b|}b\Delta(ac)\\
&-&(\Delta a)bc-(-1)^{|a|}a(\Delta b)c-(-1)^{|a|+|b|}ab(\Delta c).
\end{eqnarray*}
Together, the loop product and the BV-operator make $\H_\ast(LM)$ into a \emph{Batalin-Vilkovisky algebra}, or \emph{BV-algebra}, that we will refer to as the \emph{string topology BV-algebra of $M$}.

The string topology BV-algebra has by now been computed for several classes of manifolds.  These are the complex Stiefel manifolds \cite{\TamanoiBV} and  spheres \cite{\Menichi}, where coefficients were taken in the integers, and the complex and quaternionic projective spaces \cite{\Yang} and surfaces of genus $g>1$ \cite{\Vaintrob}, where coefficients were taken in the rational numbers.

The purpose of this paper is to compute the string topology BV-algebra, with coefficients in an arbitrary commutative ring $R$, of any compact lie group $G$.  Our main result, Theorem~\ref{MainTheorem} below, describes $\H_\ast(LG)$ in terms of four rather simpler invariants of $G$.  These are the following:
\begin{enumerate}
\item \emph{The intersection ring $\H_\ast(G)$.}\label{FirstQuantity}
The intersection ring $\H_\ast(M)$ of a closed oriented manifold $M$ is the graded-commutative ring obtained by equipping the regraded homology groups $H_{\ast+\dim M}(M)\cong H^{-\ast}(M)$ with the product Poincar\'e dual to the cup-product.

\item \emph{The action $H_\ast(G)\otimes\H_\ast(G)\to \H_\ast(G)$.}\label{SecondQuantity}  The group multiplication $\mu\colon G\times G\to G$ determines a homomorphism $H_\ast(G)\otimes H_\ast(G)\to H_\ast(G)$ that after regrading becomes $H_\ast(G)\otimes \H_\ast(G)\to \H_\ast(G)$. Thus, for $a\in H_\ast(G)$ and $x\in\H_\ast(G)$ we may form the product $ax\in\H_\ast(G)$.

\item \emph{The Hopf algebra $H_\ast(\Omega G)$.}\label{ThirdQuantity}
For any pointed space $X$, the homology groups $H_\ast(\Omega X)$ form a graded ring under the product given by concatenating loops.  When $X=G$, a theorem of Bott \cite{\EvenBott} states that these homology groups are free and concentrated in even degrees.  We may therefore equip $H_\ast(\Omega G)$ with the coproduct
\[H_\ast(\Omega G)\xrightarrow{D_\ast} H_\ast(\Omega G\times\Omega G)\cong H_\ast(\Omega G)\otimes H_\ast(\Omega G)\]
obtained using the diagonal map and the K\"unneth isomorphism.  This makes $H_\ast(\Omega G)$ into a commutative, cocommutative Hopf algebra.  We denote the coproduct of $a\in H_\ast(\Omega G)$ by $D_\ast a=\sum a_{(1)}\otimes a_{(2)}$.

\item \emph{The homology suspension $\sigma\colon H_\ast(\Omega G)\to H_{\ast+1}(G)$.}\label{FourthQuantity}
Let $X$ be a pointed space and let $\sigma\colon S^1\times\Omega X\to X$ denote the evaluation map.  The homology suspension, which we also denote by $\sigma$, is the homomorphism $H_\ast(\Omega X)\to H_{\ast+1}(X)$ defined by $\sigma(a)=\sigma_\ast([S^1]\times a)$ for $a\in H_\ast(\Omega X)$.  
\end{enumerate}
With this notation established we can state our main result.

\begin{theorem}\label{MainTheorem}
For any compact Lie group $G$ there is an isomorphism of graded rings
\begin{equation*}
\H_\ast(LG) \cong H_\ast(\Omega G) \otimes \H_\ast(G)
\end{equation*}
with respect to which the BV-operator $\Delta$ is given by
\begin{equation}\label{DeltaEquation}
\Delta(a\otimes x) = \sum a_{(1)} \otimes \sigma(a_{(2)})x.
\end{equation}
\end{theorem}

Theorem~\ref{MainTheorem} arises from the simple fact that the free loop space $LG$ splits as the product $\Omega G\times G$.  The content of the theorem lies in describing how this splitting  interacts with the string topology operations.  
We would like to note that the ring isomorphism of Theorem~\ref{MainTheorem} was stated by Tamanoi in \cite{\TamanoiBV} but not proved there.  Tamanoi also showed that $\Delta$ is a derivation of the ring $H_\ast(LG)$ whose product arises from the pointwise multiplication of loops in $G$, and then used this fact to compute the string topology BV-algebra of the special unitary groups and Stiefel manifolds.  The computation relied on a splitting of rings $H_\ast(LSU(n))\cong H_\ast(\Omega SU(n))\otimes H_\ast(SU(n))$ that has no analogue for general Lie groups.

Theorem~\ref{MainTheorem} reduces the task of computing the string topology BV-algebra of a compact Lie group $G$ to the task of computing the quantities \ref{FirstQuantity} to \ref{FourthQuantity} listed above.  In many cases these invariants of $G$ are well-known.  In particular, Bott has given a complete method for computing the third quantity, the Hopf algebra $H_\ast(\Omega G)$, in terms of the homology and cohomology of a homogeneous space called the {generating variety} \cite{\Bott}.  This, combined with Theorem~\ref{MainTheorem}, reduces the computation of the string topology BV-algebra of a Lie group to the computation of the homology and cohomology of certain finite-dimensional manifolds, and of the effect in homology and cohomology of certain maps between them.

We will see that Theorem~\ref{MainTheorem} is sufficient to give a simple and direct calculation of the string topology of the manifolds $S^1$, $S^3$ and $\mathbb{R}\mathrm{P}^3$, which are all Lie groups.  Before attempting computations for more general Lie groups, however, it is useful to explore the structure of the answer provided by Theorem~\ref{MainTheorem} in more detail.  Recall that our homology groups are taken with coefficients in a commutative ring $R$.

\begin{definition}\label{DerivationsDefinition}
Suppose either that $R$ is a field, or that $R=\Z$ and $H_\ast(G)$ is torsion-free, so that $H_\ast(G)$ becomes a Hopf algebra.  We may take a basis $p_1,\ldots,p_n$ of the odd-degree part of the primitive subspace of $H_\ast(G)$, and elements $p^1,\ldots,p^n$ of $H^\ast(G)$ for which $\langle p^i,p_j\rangle = \delta_{ij}$.   Define operators 
\[\partial_i\colon H_\ast(\Omega G)\to H_{\ast-|p_i|-1}(\Omega G), \qquad \partial_i(a) =\sum\left\langle p^i,\sigma(a_{(1)})\right\rangle a_{(2)}\]
and 
\[\delta_i\colon  \H_\ast(G)\to\H_{\ast+|p_i|}(G),\qquad \delta_i(x) = p_i x.\]
\end{definition}

\begin{proposition}\label{DerivationsProposition}
In the situation of Definition~\ref{DerivationsDefinition}, the operators $\partial_i$ and $\delta_i$ are derivations, and with respect to the isomorphism of Theorem~\ref{MainTheorem} the BV-operator $\Delta$ is given by
\begin{equation*}
\Delta = \sum_{i=1}^n \partial_i \otimes \delta_i.
\end{equation*}
\end{proposition}

This proposition makes it clear exactly to what extent $\Delta$ fails to be a derivation.  The result also makes it relatively simple to describe $\Delta$ by describing the effect of the $\partial_i$ on a set of generators of $H_\ast(\Omega G)$ and the effect of the $\delta_i$ on a set of generators of $\H_\ast(G)$.  In general, in order to describe the action of $\Delta$ on $\H_\ast(LM)$ one must describe its effect on a set of generators and on all products of pairs of these generators.  Our main application of Theorem~\ref{MainTheorem} and Proposition~\ref{DerivationsProposition} is to compute the string topology BV-algebras of the special orthogonal groups $SO(n)$ with coefficients in the rational numbers and in the integers modulo $2$.  The results are as follows.

\begin{theorem}\label{RationalSOnTheorem}
Let $m\geqslant 1$.  Then
\[\H_\ast(LSO(2m+1);\Q)\cong\Q[\alpha_0,\ldots,\alpha_{2m-1}]/U_m \otimes \Lambda_\Q[\beta_3,\beta_7,\ldots,\beta_{4m-1}]\]
and
\[\H_\ast(LSO(2m+2);\Q)\cong\Q[\alpha_0,\ldots,\alpha_{2m-1},\varepsilon_m]/U_m \otimes \Lambda_\Q[\beta_3,\beta_7,\ldots,\beta_{4m-1},\gamma_{2m+1}]\]
as graded rings.  The degrees of the generators are
\[|\alpha_i|=2i,\quad|\varepsilon_m|=2m,\quad |\beta_{4i-1}|=-4i+1,\quad |\gamma_{2m+1}|=-2m-1,\]
and in both cases the relations are given by
\[U_m=\langle \alpha_0^2-1\rangle + \langle \alpha_i^2-2\alpha_{i-1}\alpha_{i+1}+\cdots\pm 2\alpha_0\alpha_{2i}\mid 1\leqslant i\leqslant m-1\rangle.\]
With respect to the above isomorphisms, the BV-operators  on $\H_\ast(LSO(2m+1);\Q)$ and $\H_\ast(LSO(2m+2);\Q)$ are given by
\[ \Delta=\sum_{i=1}^{m}\partial_{i}\otimes\delta_{i}\qquad and \qquad \Delta=\sum_{i=1}^{m}\partial_{i}\otimes\delta_{i} + \partial_\varepsilon\otimes\delta_\varepsilon\]
respectively.  Here $\partial_{i}$, $\delta_{i}$, $\partial_\varepsilon$ and $\delta_\varepsilon$ are the derivations defined as follows:
\begin{itemize}
\item $\partial_{i}$ sends $\alpha_{j}$ to $\alpha_{j-2i+1}$ if $j\geqslant 2i-1$ and sends all other generators to $0$.
\item $\delta_{i}$ sends $\beta_{4i-1}$ to the unit and sends all other generators to 0.
\item $\partial_\varepsilon$ sends $\varepsilon_m$ to the unit and sends all other generators to 0.
\item $\delta_\varepsilon$ sends $\gamma_{2m+1}$ to the unit and sends all other generators to 0.
\end{itemize}
\end{theorem}

\begin{theorem}\label{ModTwoSOnTheorem}
Let $m\geqslant 1$.  Then as graded rings, $\H_\ast(LSO(2m+1);\Z_2)$ and $\H_\ast(LSO(2m+2);\Z_2)$ are isomorphic to
\[\Z_2[a_0,\ldots,a_{m-1},b_0,\ldots,b_{m-1}]/R_m \otimes \Z_2[c_1,c_3,\ldots,c_{2m-1}]/P_m\]
and
\[\Z_2[a_0,\ldots,a_{m-1},a_m,b_0,\ldots,b_{m-1}]/R_m \otimes \Z_2[c_1,c_3,\ldots,c_{2m+1}]/Q_m\]
respectively.  The degrees of the generators are
\[|a_i|=2i, \quad |b_i|=2m+2i, \quad |c_{2i-1}|=-2i+1,\] 
and the relations are given by
\begin{eqnarray*}
R_m&=&\langle a_0^2+1\rangle+\langle a_i^2\mid 2i\leqslant m-1\rangle+\langle a_i^2+b_{2i-m}a_0+\cdots+b_0a_{2i-m}\mid 2i\geqslant m\rangle\\
P_m&=&\langle c_{2i-1}^{r_i}\mid 2i-1\leqslant 2m\rangle\\
Q_m&=&\langle c_{2i-1}^{s_i}\mid 2i-1\leqslant 2m+1\rangle
\end{eqnarray*}
where $r_i$ is the smallest power of $2$ for which $(2i-1)r_i\geqslant(2m+1)$ and $s_i$ is the smallest power of $2$ for which $(2i-1)s_i\geqslant(2m+2)$.  With respect to the above isomorphisms the BV-operators on $\H_\ast(LSO(2m+1);\Z_2)$ and $\H_\ast(LSO(2m+2);\Z_2)$ are
\[\Delta = \sum_{i=1}^m \partial_i \otimes \delta_{i},\qquad \Delta = \sum_{i=1}^{m+1} \partial_i \otimes \delta_{i},\]
respectively, where $\partial_i$ and $\delta_{i}$ are the derivations defined as follows:
\begin{itemize}
\item $\partial_i$ sends $a_{j}$ to $a_{j-i+1}$ if $j\geqslant i-1$ and to $0$ otherwise, and similarly for the $b_j$.
\item $\delta_{i}$ sends $c_{2i-1}$ to the unit and sends all other generators to $0$.
\end{itemize}
\end{theorem}

We would like to mention two extensions of Theorem~\ref{MainTheorem}.  The first extension is to the situation where the Lie group $G$ acts smoothly on a manifold $M$.  Then the ring $H_\ast(\Omega G)$ acts on the homology groups $\H_\ast(LM)$.  It is possible to prove a result describing how this action interacts with the BV-structure on $\H_\ast(LM)$ that implies Theorem~\ref{MainTheorem} in the case $M=G$.  For the second extension, recall from \cite{\CohenJones} that if $h_\ast$ is a multiplicative homology theory in which $M$ is oriented then, just as for ordinary homology, the groups $h_{\ast+\dim M}(LM)$ admit the structure of a BV-algebra.  Theorem~\ref{MainTheorem} can be extended to give a description of $h_{\ast+\dim G}(LG)$ so long as the coefficients $h_\ast(\ast)$ are concentrated in even degrees; for more general $h_\ast$ we do not know to what extent such a result holds.

The results presented here raise the following question.  Godin has shown that the BV-structure on $\H_\ast(LM)$ can be extended to a \emph{degree $\dim M$ open-closed homological conformal field theory (HCFT) with positive boundary on the pair $(H_\ast(LM),H_\ast(M))$} \cite{\Godin}.  In particular, this endows $\H_\ast(LM)$ with a host of new operations arising from families of Riemann surfaces with boundary.  When $M=G$ it is interesting to ask whether these operations can be described in terms of the quantities \ref{FirstQuantity} to \ref{FourthQuantity} listed above and, if not, what new invariants of $G$ must be used to give a complete description.

The paper is arranged as follows.  Section~\ref{SuspensionSection} recalls some simple properties of the homology suspension that will be used throughout the rest of the paper.  In section~\ref{ExamplesSection} we illustrate Theorem~\ref{MainTheorem} by using it to compute the string topology of $S^1$, $S^3$ and $\mathbb{R}\mathrm{P}^3$ with coefficients in $\Z$.  The first two of these results are due to Menichi~\cite{\Menichi} and Tamanoi~\cite{\TamanoiBV}, but to our knowledge the third result is new.  In section~\ref{TheoremProofSection} we recall the definition of the intersection product and the string topology operations and then use these to prove Theorem~\ref{MainTheorem}, and in section~\ref{PropositionProofSection} we prove Proposition~\ref{DerivationsProposition}.  Section \ref{SOnSection} gives the proofs of Theorems~\ref{RationalSOnTheorem} and \ref{ModTwoSOnTheorem}.

\section{The homology suspension}\label{SuspensionSection}

The purpose of this short section is to recall some properties of the homology suspension that will be useful in the remainder of the paper.  The three lemmas that follow are either obvious or well-known, and we will provide references at the end of the section.  We work with homology with coefficients in a commutative ring $R$.

\begin{lemma}\label{SuspensionNaturalityLemma}
The homology suspension $\sigma\colon H_\ast(\Omega X)\to H_{\ast+1}(X)$ is natural with respect to maps of $X$.  That is, if $f\colon X\to Y$ is a continuous map of based spaces, then the diagram
\[\xymatrix{
H_\ast(\Omega X) \ar[r]^{\Omega f_\ast}\ar[d]_{\sigma} & H_\ast(\Omega Y)\ar[d]^{\sigma}\\
H_{\ast+1}(X)\ar[r]_{f_\ast} & H_{\ast+1}(Y)
}\]
commutes.
\end{lemma}

\begin{lemma}\label{SuspensionTransgressionLemma}
The homology suspension is equal to the composite
\[H_\ast(\Omega X)\to H_\ast(\Omega X,\ast)\xrightarrow{\partial_\ast^{-1}}H_{\ast+1}(PX,\Omega X)\xrightarrow{\pi_\ast}H_{\ast+1}(X,\ast)=H_{\ast+1}(X)\]
where $\Omega X\hookrightarrow PX\xrightarrow{\pi}X$ is the path-loop fibration of $X$.  In particular the homology suspension of any class $a\in H_\ast(\Omega X)$ is transgressive and its transgression is the coset $a+\partial_\ast\ker\pi_\ast$.
\end{lemma}

\begin{lemma}\label{SuspensionTheoremLemma}
Let $a,b\in H_\ast(\Omega X)$ and let $\varepsilon\colon H_\ast(\Omega X)\to R$ denote the augmentation.  Then
\[\sigma(ab) = \sigma(a)\varepsilon(b) + \varepsilon(a)\sigma(b),\]
where the product $ab$ is formed using the Pontrjagin product on $H_\ast(\Omega X)$. Also
\[D_\ast\sigma(a)=\sigma(a)\times 1 + 1\times\sigma(a),\]
where $D\colon X\to X\times X$ denotes the diagonal map.
\end{lemma}

Lemma~\ref{SuspensionNaturalityLemma} is an immediate consequence of the definitions.  Lemma~\ref{SuspensionTransgressionLemma} is an instance of the commutative diagram that follows Lemma 6.11 of \cite{\SpectralSequences}.  Lemma~\ref{SuspensionTheoremLemma} is a consequence of the Homology Suspension Theorem of \cite[Chapter VIII]{\Elements}, in whose notation the two claims are $\sigma_\ast\tau_\ast=0$ and $d_2\sigma_\ast=0$ respectively.

\section{Some simple examples}\label{ExamplesSection}

The string topology BV-algebra $\H_\ast(LS^1;\Z)$ was computed by Menichi in \cite{\Menichi}, and the string topology BV-algebra $\H_\ast(LS^3;\Z)$ was computed by Tamanoi in \cite{\TamanoiBV} and Menichi in \cite{\Menichi}.  In this section we use Theorem~\ref{MainTheorem} to give new proofs of these two results, and to compute $\H_\ast(L\mathbb{R}\mathrm{P}^3;\Z)$.  We believe that this last result is new.

\begin{proposition}[\cite{\Menichi}]
There is an isomorphism of rings
\[\H_\ast(LS^1;\Z)\cong \Z[x,x^{-1}]\otimes\Lambda_\Z[a],\qquad |x|=0,\  |a|=-1\]
under which the BV-operator $\Delta$ is given by 
\[\Delta(x^i\otimes a)=ix^i\otimes 1,\qquad \Delta(x^i\otimes 1)=0\]
for $i\in\Z$.
\end{proposition}

\begin{proof}
Since $S^1$ is the Lie group of real numbers modulo the integers, Theorem~\ref{MainTheorem} can be applied to compute the string topology BV-algebra $\H_\ast(LS^1;\Z)$.  In the rest of this proof integer coefficients should be understood.

The homology and cohomology rings of $S^1$ are $H_\ast(S^1)=\Lambda_\Z[S^1]$ and $H^\ast(S^1)=\Lambda_\Z[S^1]^\ast$.  Thus $\H_\ast(S^1)=\Lambda_\Z[a]$, where $a\in \H_{-1}(S^1)=H_0(S^1)$ denotes the class of a point and $1\in\H_0(S^1)=H_1(S^1)$ is the fundamental class $[S^1]$.  The action of $H_\ast(S^1)$ on $\H_\ast(S^1)$ is given by $[S^1]a=1$, $[S^1]1=0$.

There is a homotopy equivalence of H-spaces $\Omega S^1\simeq\Z$, with $\Id_{S^1}\in\Omega S^1$ corresponding to $1\in\Z$.  Writing the homology class of this point as $x\in H_0(\Omega S^1)$ we find that $H_\ast(\Omega S^1)=\Z[x,x^{-1}]$, that $D_\ast x = x\otimes x$, and that $\sigma(x^i)=i[S^1]$ since $x^i$ is the class of the $i$-fold map $S^1\to S^1$.

We now apply Theorem~\ref{MainTheorem}.  From the last two paragraphs we have $\H_\ast(LS^1)\cong H_\ast(\Omega S^1)\otimes\H_\ast(S^1)=\Z[x,x^{-1}]\otimes\Lambda_\Z[a]$, and from the last paragraph we have $\Delta(x^i\otimes\alpha)=x^i\otimes\sigma(x^i)\alpha=i x^i\otimes [S^1]\alpha$ for any $\alpha\in\H_\ast(S^1)$.  The stated description of $\Delta$ now follows from our description of the action of $H_\ast(S^1)$ on $\H_\ast(S^1)$.  This completes the proof.
\end{proof}

\begin{proposition}[\cite{\TamanoiBV}, \cite{\Menichi}]
There is an isomorphism of rings
\[\H_\ast(LS^3;\Z)\cong \Z[u]\otimes\Lambda_\Z[a],\qquad |u|=2,\ |a|=-3\]
under which the BV-operator $\Delta$ is given by 
\[\Delta(u^i\otimes a)=iu^{i-1}\otimes 1,\qquad \Delta(u^i\otimes 1)=0\]
for $i\geqslant 0$.
\end{proposition}
\begin{proof}
Since $S^3$ is diffeomorphic to the special unitary group $SU(2)$, Theorem~\ref{MainTheorem} can be applied to compute the string topology BV-algebra $\H_\ast(LS^3;\Z)$. Throughout the rest of the proof integer coefficients should be understood.

The homology and cohomology rings of $S^3$ are $H_\ast(S^3) =\Lambda_\Z[S^3]$ and $H^\ast(S^3) = \Lambda_\Z[S^3]^\ast$.  It follows that $\H_\ast(S^3)=\Lambda_\Z[a]$, where $a\in\H_{-3}(S^3)=H_0(S^3)$ is the homology class of a point and $1\in\H_0(S^3)=H_3(S^3)$ is the fundamental class $[S^3]$.  The action of $H_\ast(S^3)$ on $\H_\ast(S^3)$ is given by $[S^3]a=1$, $[S^3]1=0$.

A simple argument using the Serre spectral sequence of the path-loop fibration $\Omega S^3\to P S^3\to S^3$ shows that $H_\ast(\Omega S^3)$ is isomorphic to the polynomial ring $\Z[u]$, where $u\in H_2(\Omega S^3)$ is the transgression of $[S^3]\in H_3(S^3)$.  It follows from Lemma~\ref{SuspensionTransgressionLemma} that $\sigma(u)=[S^3]$, and for degree reasons that $\sigma(u^i)=0$ for $i\neq 1$.  Also for degree reasons we have $D_\ast u=u\otimes 1 + 1\otimes u$, so that $D_\ast u^i=\sum\binom{i}{j}u^{i-j}\otimes u^j$ and consequently $\sum {u^i}_{(1)}\otimes \sigma({u^i}_{(2)})=iu^{i-1}\otimes[S^3]$ for all $i\geqslant 0$.

We now apply Theorem~\ref{MainTheorem}.  The isomorphism $\H_\ast(LS^3)\cong H_\ast(\Omega S^3)\otimes\H_\ast(S^3)\cong\Z[u]\otimes\Lambda_\Z[a]$ follows from the last two paragraphs.  From the last paragraph we have $\Delta(u^i\otimes\alpha)=\sum {u^i}_{(1)}\otimes \sigma({u^i}_{(2)})\alpha=iu^{i-1}\otimes [S^3]\alpha$.  The stated description of $\Delta$ now follows from our description of the action of $H_\ast (S^3)$ on $\H_\ast(S^3)$.
\end{proof}

\begin{proposition}
There is an isomorphism of rings
\[\H_\ast(L\mathbb{R}\mathrm{P}^3;\Z)\cong \Z[u,v]\otimes\Lambda_\Z[a,b] / \langle v^2-1,2b,ab\rangle\]
where the degrees of the generators are
\[|u|=2,\quad |v|=0,\quad |a|=-3,\quad |b|=-2.\]
Under this isomorphism the BV-operator $\Delta$ is given by
\begin{gather*}
\Delta(u^iv^j\otimes a) = 2iu^{i-1}v^j\otimes 1+ ju^iv^j\otimes b,\quad
\Delta(u^iv^j\otimes b)=0,\quad\Delta(u^iv^j\otimes 1)=0.
\end{gather*}
\end{proposition}

\begin{proof}
Since $\mathbb{R}\mathrm{P}^3$ is diffeomorphic to the quotient of $S^3\cong{SU}(2)$ by its centre, Theorem~\ref{MainTheorem} can be applied to compute the string topology BV-algebra $\H_\ast(L\mathbb{R}\mathrm{P}^3;\Z)$.  Throughout this proof integer coefficients should be understood.

The homology and cohomology rings of $\mathbb{R}\mathrm{P}^3$ are given by
\begin{align*}
H_\ast(\mathbb{R}\mathrm{P}^3)&= \Lambda_\Z\left[\rho,[\mathbb{R}\mathrm{P}^3]\right]/\left\langle 2\rho,\rho[\mathbb{R}\mathrm{P}^3] \right\rangle,\\
H^\ast(\mathbb{R}\mathrm{P}^3)&= \Lambda_\Z\left[\tau,[\mathbb{R}\mathrm{P}^3]^\ast\right]/\left\langle 2\tau,\tau[\mathbb{R}\mathrm{P}^3]^\ast \right\rangle,
\end{align*}
where $\rho$ is a generator of $H_1(\mathbb{R}\mathrm{P}^3)$ and $\tau$ is a generator of $H^2(\mathbb{R}\mathrm{P}^3)$.  Consequently $\H_\ast(\mathbb{R}\mathrm{P}^3) = \Lambda_\Z[a,b]/\langle 2b,ab\rangle$, where $a\in \H_{-3}(\mathbb{R}\mathrm{P}^3)=H_0(\mathbb{R}\mathrm{P}^3)$ is the homology class of a point and $b\in\H_{-2}(\mathbb{R}\mathrm{P}^3)=H_1(\mathbb{R}\mathrm{P}^3)$ is a generator.  The action of $H_\ast(\mathbb{R}\mathrm{P}^3)$ on $\H_\ast(\mathbb{R}\mathrm{P}^3)$ is given by
\[[\mathbb{R}\mathrm{P}^3]a=1,\quad [\mathbb{R}\mathrm{P}^3]b=0,\quad [\mathbb{R}\mathrm{P}^3]1=0,\quad \rho a =b,\quad\rho b=0,\quad\rho 1=0.\]

Since $\mathbb{R}\mathrm{P}^3$ is the quotient $S^3/\Z_2$, there is a homotopy equivalence of H-spaces $\Omega \mathbb{R}\mathrm{P}^3\simeq\Omega S^3\times\Z_2$ and a corresponding ring isomorphism $H_\ast(\Omega\mathbb{R}\mathrm{P}^3)\cong\Z[u]\otimes\Z[v]/\langle v^2-1\rangle$.  If we write $p\colon S^3\to\mathbb{R}\mathrm{P}^3$ for the quotient map then $u\in H_2(\Omega \mathbb{R}\mathrm{P}^3)$ is equal to $\Omega p_\ast u$, where by abuse of notation $u\in H_2(\Omega S^3)$ is the class described in the proof of the proposition above.  By Lemma~\ref{SuspensionNaturalityLemma} we have $\sigma(u)=\sigma(\Omega p_\ast u)=p_\ast\sigma(u)=p_\ast[S^3]=2[\mathbb{R}\mathrm{P}^3]$ and by naturality of the diagonal we have $D_\ast u =u\otimes 1+1\otimes u$.  The class $v\in H_0(\Omega\mathbb{R}\mathrm{P}^3)$ is the homology class of point in $\Omega\mathbb{R}\mathrm{P}^3$ corresponding to a noncontractible loop in $\mathbb{R}\mathrm{P}^3$ and so has $\sigma(v)=\rho$ and $D_\ast v=v\otimes v$.  It now follows from Lemma~\ref{SuspensionTheoremLemma} that $\sigma(v^j)=j\rho$, that $\sigma(uv^j)=2[\mathbb{R}\mathrm{P}^3]$, and that $\sigma(u^iv^j)=0$ for $i\geqslant 2$.  We therefore have
\begin{eqnarray*}
\sum (u^iv^j)_{(1)}\otimes\sigma((u^iv^j)_{(2)})
&=&\sum{\textstyle\binom{i}{j}u^{i-k}v^j}\otimes \sigma(u^kv^j)\\
&=&iu^{i-1}v^j\otimes2[\mathbb{R}\mathrm{P}^3]+ju^iv^j\otimes\rho.
\end{eqnarray*}

We now apply Theorem~\ref{MainTheorem}.  The isomorphism $\H_\ast(L\mathbb{R}\mathrm{P}^3)\cong H_\ast(\Omega \mathbb{R}\mathrm{P}^3)\otimes\H_\ast(\mathbb{R}\mathrm{P}^3)=\Z[u,v]\otimes\Lambda_\Z[a,b] / \langle v^2-1,2b,ab\rangle$ follows from the last two paragraphs, and from the last paragraph we have $\Delta(u^iv^j\otimes\alpha)=iu^{i-1}v^j\otimes 2[\mathbb{R}\mathrm{P}^3]\alpha + ju^iv^j\otimes\rho\alpha$ for any $\alpha\in\H_\ast(\mathbb{R}\mathrm{P}^3)$.  The stated description of $\Delta$ now follows from our description of the action of $H_\ast(\mathbb{R}\mathrm{P}^3)$ on $\H_\ast(\mathbb{R}\mathrm{P}^3)$.
\end{proof}

\section{Proof of Theorem~\ref{MainTheorem}}\label{TheoremProofSection}
We continue to work with homology with coefficients in a commutative ring $R$.
There is a homeomorphism
\[\Theta\colon\Omega G\times G\xrightarrow{\cong} LG\]
defined by $\Theta(\delta,g)(t)=\delta(t)g$ for $\delta\in\Omega G$, $g\in G$ and $t\in S^1$.  The theorem of Bott \cite{\EvenBott} which we mentioned in the introduction states that $H_\ast(\Omega G)$ is a free $R$-module concentrated in even degrees.  It follows that the Kunneth isomorphism
\[H_\ast(\Omega G)\otimes H_\ast(G)\xrightarrow{\cong}H_\ast(\Omega G\times G)\]
holds.  Combining these two facts and applying the degree-shift we have an isomorphism of $R$-modules
\[\Phi\colon H_\ast(\Omega G)\otimes \H_\ast(G)\xrightarrow{\cong}\H_\ast(LG)\]
defined by $\Phi(a\otimes x)=\Theta_\ast(a\times x)$ for $a\in H_\ast(\Omega G)$ and $x\in H_\ast(G)$.  In this section we will prove Theorem~\ref{MainTheorem} by showing that $\Phi$ is a ring-homomorphism, and then showing that after applying $\Phi$ the BV-operator $\Delta$ is given by equation \eqref{DeltaEquation}.

\subsection{The intersection product and the string topology operations.}\label{DefinitionsSubsection}

We begin by recalling the construction of the intersection product on $\H_\ast(M)$ when $M$ is a closed oriented manifold of dimension $m$.  We will use the construction described by Cohen and Jones \cite{\CohenJones}.  Let $D\colon M\to M\times M$ denote the diagonal map.  Since $D$ is an embedding of manifolds with normal bundle $TM$, there is a tubular neighbourhood $\nu_D\subset M\times M$ of $D(M)$ diffeomorphic to the total space $TM$.  There is an associated Pontrjagin-Thom collapse map $D_!\colon M\times M\to M^{TM}$
and a map in homology that we denote in the same way:
\[D_!\colon H_\ast(M\times M)\to H_\ast (M^{TM}).\]
We also have the Thom isomorphism
\[\Th\colon\tilde H_\ast(M^{TM}){\cong} H_{\ast-m}(M).\]
With this notation established, the intersection product of $x,y\in\H_\ast(M)=H_{\ast+m}(M)$ is given by
\begin{equation}\label{IntersectionEquation}x\cdot y = (-1)^{m|y|+m}\Th\circ D_!(x\times y).\end{equation}
Note that here the symbol $|y|$ denotes the degree of $y$ as an element of $\H_\ast(M)$.

We now recall the construction of the loop product on $\H_\ast(LM)$.   We again use the construction due to Cohen and Jones \cite{\CohenJones}.  Recall that the free loop space $LM$ is the total space of a fibre bundle $\ev\colon LM\to M$ that sends a loop in $M$ to its value at the basepoint $0\in S^1$.  Write $L^2M=\{(\delta_1,\delta_2)\mid\ev(\delta_1)=\ev(\delta_2)\}$ for the space of pairs of composable loops in $M$, $\ev\colon L^2M\to M$ for the map that sends such a pair to their common basepoint, $\tilde D\colon L^2M\to LM\times LM$ for the inclusion, and $\gamma\colon L^2M\to LM$ for the map that composes loops, so that for $(\delta_1,\delta_2)\in L^2M$ we have 
\[\gamma(\delta_1,\delta_2)(t)=\left\{\begin{array}{lc}\delta_1(2t), & 0\leqslant t\leqslant 1/2,\\ \delta_2(2t-1), &1/2\leqslant t\leqslant 1.\end{array}\right.\]
The diagram
\[\xymatrix{
L^2M \ar[r]^-{\tilde D}\ar[d]_\ev & LM\times LM \ar[d]^{\ev\times\ev}\\
M\ar[r]_-{D} & M\times M
}\]
is a pullback square whose vertical maps are the projections of fibre bundles.  As described in the last paragraph, $D$ is an embedding of manifolds with normal bundle $TM$ and there is a tubular neighbourhood $\nu_D\subset M\times M$ of $D(M)$ diffeomorphic to the total space $TM$.  It follows that $\tilde D(L^2M)$ admits a tubular neighbourhood $\nu_{\tilde D}$ homeomorphic to $\ev^\ast TM$.  There is an associated Pontrjagin-Thom collapse map $\tilde D_!\colon LM\times LM\to L^2M^{\ev^\ast TM}$  and a map in homology that we denote in the same way:
\[\tilde D_!\colon H_\ast(LM\times LM)\to H_\ast(L^2M^{\ev^\ast TM}).\]
We also have the Thom isomorphism
\[\Th\colon\tilde H_\ast(L^2M^{\ev^\ast TM}){\cong} H_{\ast-m}(L^2M).\]
With this notation established, the loop product of $x,y\in\H_\ast(LM)=H_{\ast +m}(LM)$ is given by
\begin{equation}\label{LoopEquation}x\cdot y = (-1)^{m|y|+m}\gamma_\ast\circ\Th\circ \tilde D_!(x\times y).\end{equation}
Note that here the symbol $|y|$ denotes the degree of $y$ as an element of $\H_\ast(LM)$.

Finally we recall the definition of the BV-operator $\Delta$.  The circle $S^1$ acts on $LM$ by rotating loops, or in other words by the action $\rho\colon S^1\times LM\to LM$ defined by $\rho(s,\delta)(t)=\delta(s+t)$ for $\delta\in LM$ and $s,t\in S^1$.  The BV-operator
\[\Delta\colon\H_\ast(LM)\to\H_{\ast+1}(LM)\]
is defined by $\Delta(x)=\rho_\ast([S^1]\times x)$ for $x\in\H_\ast(LM)$.

\begin{remark}
The sign correction $(-1)^{m|y|+m}$ appearing in \eqref{IntersectionEquation} and \eqref{LoopEquation} is necessary if one wishes to obtain a graded associative, graded commutative product.  For example, it is routine to verify that $\Th\circ D_!(y\times x)=(-1)^{m|x|+m|y|+|x||y|}\Th\circ D_!(x\times y)$, so that the product defined in \eqref{IntersectionEquation} satisfies $x\cdot y=(-1)^{|x||y|}y\cdot x$.  See the discussion in \cite[\S 4.6]{\Godin}.

In fact Cohen and Jones \cite{\CohenJones} use the collapse maps $D_!$, $\tilde D_!$ to describe ring structures on the Thom spectra $M^{-TM}$, $LM^{-\ev^\ast TM}$.  These lead to ring structures on $H_\ast(M^{-TM})$, $H_\ast(LM^{-\ev^\ast TM})$ that after applying the Thom isomorphisms
\[H_\ast(M^{-TM})\cong\H_\ast(M),\qquad H_\ast(LM^{-\ev^\ast TM})\cong\H_\ast (LM)\]
give ring structures on $\H_\ast(M)$ and $\H_\ast(LM)$.
It is simple to show that these products on $\H_\ast(M)$, $\H_\ast(LM)$ are given by \eqref{IntersectionEquation} and \eqref{LoopEquation}.
\end{remark}

\subsection{Proof that $\Phi$ is a ring-isomorphism.}
By the definition of the module isomorphism $\Phi$ and the formulas \eqref{IntersectionEquation}, \eqref{LoopEquation} of the last subsection, and recalling that $H_\ast(\Omega G)$ is concentrated in even degrees, the claim that $\Phi$ is a ring isomorphism is equivalent to the claim that for $a,a'\in H_\ast(\Omega G)$ and $x,x'\in H_\ast (G)$ we have
\[\gamma_\ast\circ\Th\circ\tilde D_!\left(\Theta_\ast(a\times x)\times\Theta_\ast(a'\times x')\right)=\Theta_\ast\left(p_\ast(a\times a')\times\Th\circ D_!(x\times x')\right)\]
where $p\colon\Omega G\times\Omega G\to\Omega G$ is the concatenation of based loops.  Recall that we have the homeomorphism $\Theta\colon\Omega G\times G\to LG$.  Write $\Theta^2\colon\Omega G\times\Omega G\times G\to L^2G$ for the homeomorphism defined by $\Theta^2(\delta_1,\delta_2,g)=(\Theta(\delta_1,g),\Theta(\delta_2,g))$ for $\delta_1,\delta_2\in\Omega G$ and $g\in G$.  We will now treat the isomorphisms $\Theta$, $\Theta^2$ as identifications and work through the definition of the loop product on $\H_\ast(\Omega G\times G)$, which we recalled in the last subsection, in order to prove the equation above.

First, $\ev\colon\Omega G\times G\to G$ is just the projection from a product to one of its factors, so the Pontrjagin-Thom map $\tilde D_!\colon(\Omega G\times G)\times(\Omega G\times G)\to(\Omega G\times\Omega G\times G)^{\ev^\ast TG}$ is the composition
\begin{multline*}
(\Omega G\times G)\times(\Omega G\times G)\cong\Omega G\times\Omega G\times (G\times G)\\
\xrightarrow{\Id\times\Id\times D_!}\Omega G\times \Omega G\times G^{TG}\xrightarrow{c}(\Omega G\times\Omega G\times G)^{\ev^\ast TG},\end{multline*}
where the first map shuffles the factors and the last map $c$ collapses $\Omega G\times\Omega G\times\ast$ to a point.  Thus $\tilde D_!((a\times x)\times (a'\times x'))=c_\ast (a\times a'\times D_!(x\times x'))$.  

Second, $\ev\colon\Omega G\times\Omega G\times G\to G$ is the projection from a product to one of its factors and so the composition
\[H_\ast(\Omega G\times\Omega G\times G^{TG})\xrightarrow{c_\ast}H_\ast((\Omega G\times\Omega G\times G)^{\ev^\ast TG})\xrightarrow{\Th} H_{\ast-\dim G}(\Omega G\times \Omega G\times G)\]
is given by $\Th\circ c_\ast(a\times a'\times u)=a\times a'\times\Th(u)$ for any $u\in H_\ast (G)$.  In particular $\Th\circ c_\ast(a\times a'\times D_!(x\times x'))=a\times a'\times \Th\circ D_!(x\times x')$.

Finally, $\gamma\colon\Omega G\times\Omega G\times G\to\Omega G\times G$ is just the product $p\times\Id$.  Thus $\gamma_\ast(a\times a'\times\Th\circ D_!(x\times x'))=p_\ast(a\times a')\times\Th\circ D_!(x\times x')$.  This, together with the conclusions of the last two paragraphs, proves the claim.\qed

\subsection{Proof of equation \eqref{DeltaEquation}.}
We must prove that with respect to the ring isomorphism $\Phi$ the BV-operator $\Delta$ is given by equation \eqref{DeltaEquation}.  In light of the definition of $\Phi$ and the description of $\Delta$ given in \S\ref{DefinitionsSubsection} we must prove that the equation
\[\rho_\ast\left([S^1]\times\Theta_\ast(a\times x)\right)=\sum \Theta_\ast\left(a_{(1)}\times\sigma(a_{(2)})x\right)\]
holds for any $a\in H_\ast(\Omega G)$ and $x\in H_\ast(G)$.

Let $\tau\colon S^1\times\Omega G\to\Omega G$ be the map defined by $\tau(s,\delta)(t)=\delta(s+t)\delta(s)^{-1}$ for $\delta\in\Omega G$ and $s,t\in S^1$.
Then note that $\rho(s,\Theta(\delta,g))(t)=\Theta(\tau(s,\delta),\delta(s)g)(t)$ for $s,t\in S^1$, $\delta\in\Omega G$ and $g\in G$. In other words $\rho\circ(\Id\times\Theta)$ is the composite
\begin{multline*}
S^1\times\Omega G\times G\xrightarrow{D\times D\times\Id}(S^1\times S^1)\times(\Omega G\times\Omega G)\times G\\
\cong(S^1\times\Omega G)\times(S^1\times\Omega G)\times G\xrightarrow{\tau\times\sigma\times\Id}\Omega G\times G\times G\\
\xrightarrow{\Id\times\mu}\Omega G\times G\xrightarrow{\Theta} LG,
\end{multline*}
where the isomorphism shuffles the factors and $\mu\colon G\times G\to G$ denotes the group multiplication.  We must therefore prove that this composite sends the homology class $[S^1]\times a\times x\in H_\ast(S^1\times\Omega G\times G)$ to $\sum \Theta_\ast\left(a_{(1)}\times\sigma(a_{(2)})x\right)$.

The map $S^1\times\Omega G\times G\to(S^1\times\Omega G)\times(S^1\times \Omega G)\times G$ in the above composite sends the homology class $[S^1]\times a\times x$ to
\begin{equation}\label{HomologyClass}\sum ([S^1]\times a_{(1)})\times (1\times a_{(2)})\times x+\sum(1\times a_{(2)})\times ([S^1]\times a_{(2)})\times x.\end{equation}
Recall that $\sigma_\ast([S^1]\times a_{(2)})=\sigma(a_{(2)})$.  Note that $\tau_\ast([S^1]\times a_{(1)})=0$ because $H_\ast(\Omega G)$ is concentrated in even degrees and that $\tau_\ast(1\times a_{(2)})=a_{(2)}$.   It follows that $\tau\times\sigma\times\Id$ sends the homology class \eqref{HomologyClass} to $\sum a_{(1)}\times\sigma(a_{(2)})\times x$.  This class is, in turn, sent by $\Theta\circ(\Id\times\mu)$ into the class $\sum \Theta_\ast(a_{(1)}\times\sigma(a_{(2)})x)$.  In other words the composite above sends $[S^1]\times a\times x$ to $\sum a_{(1)}\times\sigma(a_{(2)})x$ as required.\qed

\section{Proof of Proposition~\ref{DerivationsProposition}}\label{PropositionProofSection}

In this section we will prove Proposition~\ref{DerivationsProposition}.  This is a routine consequence of the properties of the homology suspension listed in Lemma~\ref{SuspensionTheoremLemma} once we have established the following property of the intersection product, which is valid for homology with coefficients in any commutative ring $R$.  To state it, we assume that $G$ acts smoothly on a closed $m$-dimensional manifold $M$ and that this action preserves the orientation of $M$.  This action makes $\H_\ast(M)$ into a module over the ring $H_\ast(G)$.

\begin{lemma}\label{GIntersectionProductLemma}
Let $x,y\in\H_\ast(M)$ and let $\alpha\in H_\ast(G)$ be such that $D_\ast \alpha=\sum\alpha_{(1)}\times\alpha_{(2)}$.  Then
\[\alpha(x\cdot y) = \sum (-1)^{|\alpha_{(2)}||x|} (\alpha_{(1)}x)\cdot(\alpha_{(2)}y).\]
\end{lemma}
\begin{proof}
To prove this claim we will use the description of the intersection product in terms of the Pontrjagin-Thom construction, which we recalled in \S\ref{DefinitionsSubsection}.  Write $\mu_1\colon G\times M\to M$ for the action, $\mu_2\colon G\times M\times M\to M\times M$ for the diagonal action $\mu_2(g,m_1,m_2)=(\mu_1(g,m_1),\mu_2(g,m_2))$ and write $\mu_3\colon G\times M^{TM}\to M^{TM}$ for the action that preserves the point at infinity and that restricts to the action $G\times TM\to TM$ obtained by differentiating $\mu_1$.

Since $D\colon M\to M\times M$ is $G$-equivariant, the tubular neighbourhood of $D(M)$ may be chosen in a $G$-equivariant way, and it follows that $\mu_3\circ(\Id\times D_!)=D_!\circ\mu_2$.  Further, since the action of $G$ on $M$ preserves orientations, the two oriented bundles $\pi_2^\ast TM$ and $\mu_1^\ast TM$ over $G\times M$ are isomorphic, and consequently $\Th\circ{\mu_3}_\ast(a\times z)={\mu_1}_\ast (a\times\Th(z))$ for any $z\in H_\ast(M^{TM})$.  Thus
\begin{eqnarray*}
\alpha(x\cdot y)
&=&(-1)^{m|y|+m}\mu_\ast(\alpha\times\Th\circ D_!(x\times y))\\
&=&(-1)^{m|y|+m}\Th\circ{\mu_3}_\ast(\alpha\times D_!(x\times y))\\
&=&(-1)^{m|y|+m}\Th\circ D_!\circ{\mu_2}_\ast(\alpha\times x\times y)\\
&=&\sum(-1)^{m|y|+|\alpha_{(2)}|(|x|+m)+m}\Th\circ D_!(\alpha_{(1)}x\times\alpha_{(2)}y)\\
&=&\sum(-1)^{|\alpha_{(2)}||x|}(\alpha_{(1)}x)\cdot(\alpha_{(2)}y)
\end{eqnarray*}
as required.
\end{proof}

We now prove Proposition~\ref{DerivationsProposition}.  Let $a\in H_\ast(\Omega G)$ and let $x\in\H_\ast(G)$.  By Lemma~\ref{SuspensionTheoremLemma}, the classes $\sigma(a_{(2)})$ are all primitive, and so $\sigma(a_{(2)})=\sum_i\left\langle p^i,\sigma(a_{(2)})\right\rangle p_i$.  It follows from this that $\Delta(a\otimes x)= \sum_i\partial_i a\otimes\delta_i x$ as required.  
It remains  to show that each of the $\delta_i$ and $\partial_i$ is a derivation.  Since each $p_i$ is primitive and of odd degree, the fact that each $\delta_i$ is a derivation is an immediate consequence of Lemma~\ref{GIntersectionProductLemma} above.  Finally we must show that each $\partial_i$ is a derivation.  Let $a,b\in H_\ast(\Omega G)$.  We must show that $\partial_i(ab)=(\partial_i a)b+a(\partial_i b)$.  We may assume without loss that each of $a$, $b$ is concentrated in a single component of $\Omega G$, so that we may write 
\begin{align*}
D_\ast a&= a\otimes\gamma_a + \sum a_{(1)}^+\otimes a_{(2)}^+,\\
D_\ast b&= b\otimes \gamma_b + \sum b_{(1)}^+\otimes b_{(2)}^+,
\end{align*}
where $\gamma_a,\gamma_b$ are points of $\Omega G$  (and, by abuse of notation, the homology classes of those points) and each of the $a_{(2)}^+$ and  $b_{(2)}^+$ has positive degree.  Then using the fact that $D_\ast(ab)=(D_\ast a)(D_\ast b)$, and using Lemma~\ref{SuspensionTheoremLemma} to show that $\sigma(\gamma_a\gamma_b)=\sigma(\gamma_a)+\sigma(\gamma_b)$, $\sigma(\gamma_a b_{(2)}^+)=\sigma(b_{(2)}^+)$, $\sigma(a_{(2)}^+\gamma_b)=\sigma(a_{(2)}^+)$ and $\sigma(a_{(2)}^+b_{(2)}^+)=0$, it follows from the definitions that $\partial_i(ab)=(\partial_ia)b+a(\partial_ib)$.  This completes the proof.\qed

\section{String topology of special orthogonal groups}\label{SOnSection}

In this section we will prove Theorems \ref{RationalSOnTheorem} and \ref{ModTwoSOnTheorem}, which describe the string topology BV-algebra of $SO(n)$ for $n\geqslant 3$ and with coefficients in $\Q$ and $\Z_2$.  The proof will proceed using Theorem~\ref{MainTheorem} and Proposition~\ref{DerivationsProposition}.  We must therefore compute the quantities \ref{FirstQuantity} to \ref{FourthQuantity} listed in the introduction.
First, in \S\ref{HomologySubsection} we recall the homology and cohomology of $SO(n)$ with coefficients in $\Q$ and $\Z_2$ and we use this to describe quantities \ref{FirstQuantity} and \ref{SecondQuantity} from the introduction.  These results are well-known, and we will be referring to the treatment given in \cite{\Hatcher}.  Next, in \S\ref{OmegaZeroSubsection} we recall Bott's computation of $H_\ast(\Omega_0SO(n);\Z)$ from \cite{\Bott} and in \S\ref{OmegaSubsection} we use this to describe quantity \ref{ThirdQuantity} from the introduction.  Then in \S\ref{SuspensionSubsection} we compute the homology suspension using the results of the earlier subsections and some Serre spectral sequences.  Finally, in \S\ref{ProofSubsection} we prove Theorems~\ref{RationalSOnTheorem} and \ref{ModTwoSOnTheorem} by combining the results of the earlier subsections.

\subsection{Homology of $SO(n)$.}\label{HomologySubsection}
In this subsection we will take coefficients in either $\Z_2$ or $\Q$ and describe the rings $H_\ast(SO(n))$ and $\H_\ast(SO(n))$.  We will also describe the derivations of $\H_\ast(SO(n))$ determined as in Proposition~\ref{DerivationsProposition} by a basis of odd-degree primitives in $H_\ast(SO(n))$.  The results are well-known, and we refer throughout to section 3.D of \cite{\Hatcher}.  We summarize the main points as follows:

\begin{proposition}\label{RationalHomologyProposition}
Let $m\geqslant 1$.  Then
\begin{gather*}
\H_\ast(SO(2m+1);\Q)=\Lambda_\Q[\beta_3,\ldots,\beta_{4m-1}],\\
\H_\ast(SO(2m+2);\Q)=\Lambda_\Q[\beta_3,\ldots,\beta_{4m-1},\gamma_{2m+1}],
\end{gather*}
where the degrees of the generators are as in Theorem~\ref{RationalSOnTheorem}.

There is a basis $a_{3},a_7,\ldots,a_{4m-1}$ of the odd degree primitive subspace of $H_\ast(SO(2m+1);\Q)$ and a basis $a_3,a_7,\ldots,a_{4m-1},b_{2m+1}$ of the odd degree primitive subspace of $H_\ast(SO(2m+2);\Q)$.  The corresponding derivations $\delta_1,\ldots,\delta_m$ of $\H_\ast(SO(2m+1);\Q)$ and $\delta_1,\ldots,\delta_m,\delta_\varepsilon$ of $\H_\ast(SO(2m+2);\Q)$ are as described in Theorem~\ref{RationalSOnTheorem}.
\end{proposition}

\begin{proposition}\label{ModTwoHomologyProposition}
Let $m\geqslant 1$.  Then
\begin{gather*}
\H_\ast(SO(2m+1);\Z_2)=\Z_2[c_1,c_3,\ldots,c_{2m-1}]\big\slash P_{m},\\
\H_\ast(SO(2m+2);\Z_2)=\Z_2[c_1,c_3,\ldots,c_{2m+1}]\big\slash Q_m,
\end{gather*}
where the relations and the degrees of the generators are as in Theorem~\ref{ModTwoSOnTheorem}.  There is a basis $q_1,q_3,\ldots,q_{2m-1}$ of the odd primitive subspace of $H_\ast(SO(2m+1);\Z_2)$ and a basis $q_1,q_3,\ldots,q_{2m+1}$ of the odd primitive subspace of $H_\ast(SO(2m+2);\Z_2)$.
The corresponding derivations $\delta_1,\ldots,\delta_{m}$ of $\H_\ast(SO(2m+1);\Z_2)$ and $\delta_1,\ldots,\delta_{m+1}$ of $\H_\ast(SO(2m+2);\Z_2)$ are as described in Theorem~\ref{ModTwoSOnTheorem}.
\end{proposition}

We begin by recalling the rational homology and cohomology of $SO(2m+1)$ and $SO(2m+2)$ for $m\geqslant 1$. First, according to \cite{\Hatcher} the rational homology ring of $SO(2m+1)$ is given by 
\[H_\ast(SO(2m+1);\Q)=\Lambda_\Q[a_3,a_7,\ldots,a_{4m-1}]\]
where the generators are primitive and lie in degrees $|a_{4i-1}|=4i-1$. Because the generators are primitive we have
\[H^\ast(SO(2m+1);\Q)=\Lambda_\Q[a_3^\ast,\ldots,a_{4m-1}^\ast]\]
where the dual elements are formed with respect to the basis of monomials.  It follows from the construction in \cite{\Hatcher} that the fundamental class of $SO(2m+1)$ is equal to $a_3\cdots a_{4m-1}$.  Second, according to \cite{\Hatcher} the rational homology ring of $SO(2m+2)$ is given by 
\[H_\ast(SO(2m+2);\Q)=\Lambda_\Q[a_3,a_7,\ldots,a_{4m-1},b_{2m+1}]\]
where the generators are primitive and lie in degrees $|a_{4i-1}|=4i-1$, $|b_{2m+1}|=2m+1$. Because the generators are primitive we have
\[H^\ast(SO(2m+2);\Q)=\Lambda_\Q[a_3^\ast,\ldots,a_{4m-1}^\ast,b_{2m+1}^\ast]\]
where the dual elements are formed with respect to the basis of monomials.  It follows from the construction in \cite{\Hatcher} that the fundamental class of $SO(2m+2)$ is equal to $a_3\cdots a_{4m-1}b_{2m+1}$.

\begin{proof}[Proof of Proposition~\ref{RationalHomologyProposition}.]
The descriptions of the intersection rings are immediate from the above description of the cohomology rings and the fundamental class by setting
\[\beta_{4i-1}=a_{4i-1}^\ast\cap(a_3\cdots a_{4m-1})=(-1)^{i-1}a_3\cdots\widehat{a_{4i-1}}\cdots a_{4m-1}\]
in the first case and
\begin{align*}
\beta_{4i-1}&=a_{4i-1}^\ast\cap(a_3\cdots a_{4m-1}b_{2m+1})=(-1)^{i-1}a_3\cdots\widehat{a_{4i-1}}\cdots a_{4m-1}b_{2m+1},\\
\gamma_{2m+1}&=b_{2m+1}^\ast\cap(a_3\cdots a_{4m-1}b_{2m+1})=(-1)^m a_3\cdots a_{4m-1}.
\end{align*}
in the second case.  The bases of odd-degree primitives are immediate from the above description of the homology rings, as are the descriptions of the corresponding derivations.
\end{proof}

We now move on to the case of $\Z_2$ coefficients.  Recall from \cite{\Hatcher} that for $n\geqslant 2$ the homology and cohomology rings of $SO(n+1)$ are given by
\begin{align*}H_\ast(SO(n+1);\Z_2)&=\Lambda_{\Z_2}[e_1,e_2,\ldots,e_{n}], & |e_i|&=i,\\
H^\ast(SO(n+1);\Z_2)&=\Z_2[\beta_1,\beta_3,\ldots]\big\slash\langle\beta_i^{p_i}\rangle, & |\beta_{2i-1}|&=2i-1,
\end{align*}
where each $\beta_{2i-1}$ is the dual of $e_{2i-1}$ with respect to the basis of monomials and $p_i$ is the smallest power of $2$ for which $p_i(2i-1)\geqslant (n+1)$.  The $\beta_{2i-1}$ are all primitive.  It also follows from \cite{\Hatcher} that the fundamental class of $SO(n+1)$ is the product $e_1\cdots e_{n}$.

\begin{proof}[Proof of Proposition~\ref{ModTwoHomologyProposition}.]
The description of the intersection ring follows immediately from the description of the cohomology ring given above by setting
\[ c_{2i-1}=\beta_{2i-1}\cap (e_1\cdots e_n)=e_1\cdots \widehat{e_{2i-1}}\cdots e_n.\]
Since the primitive subspace of $H_\ast(SO(n+1);\Z_2)$ is dual to the quotient of $H^\ast(SO(n+1);\Z_2)$ by its decomposable elements, it follows that the primitive subspace has a basis $q_1,q_3,\ldots,q_{2m+1}$, $m=[n/2]$, uniquely determined by $\langle\beta_{2i-1},q_{2j-1}\rangle=\delta_{ij}$.  It is clear that $q_{2i-1}=e_{2i-1}$ modulo decomposables and that the product of any decomposable with any of the $c_{2i-1}$ vanishes.  From this the description of the derivation $\delta_{i}$ follows immediately.
\end{proof}

\begin{note}\label{HomologyNote}
Before ending this section we record for later use the following facts. These will be crucial for our computation of the homology suspension.  First, in rational homology:
\begin{itemize}
\item The projection $SO(2m+1)\to SO(2m+1)/SO(2m-1)$ sends the generator $a_{4m-1}$ to the fundamental class $[SO(2m+1)/SO(2m-1)]$.
\item The projection $SO(2m+2)\to S^{2m+1}$ sends the generator $b_{2m+1}$ to the fundamental class $[S^{2m+1}]$ and sends all $a_{4i-1}$ to $0$.
\item The Hopf algebra maps $H_\ast(SO(2m+1);\Q)\to H_\ast(SO(2m+2);\Q)$ and $H_\ast(SO(2m+1);\Q)\to H_\ast(SO(2m+3);\Q)$ induced by the inclusions $SO(2m+1)\hookrightarrow SO(2m+2)$ and $SO(2m+1)\hookrightarrow SO(2m+3)$ can both be described by $a_{4i-1}\mapsto a_{4i-1}$ for $i=1,\ldots, m$.  In particular, both are injections.
\end{itemize}
And in $\Z_2$ homology:
\begin{itemize}
\item The map $SO(n+1)\to S^n$ sends the generator $e_n$ to the fundamental class $[S^n]$.  In particular, $SO(2m+2)\to S^{2m+1}$ sends $q_{2m+1}$ to $[S^{2m+1}]$.
\item The Hopf algebra map $H_\ast(SO(n+1);\Z_2)\to H_\ast(SO(n+2);\Z_2)$ induced by the inclusion $SO(n+1)\hookrightarrow SO(n+2)$ is the injection given by $e_i\mapsto e_i$ for $i=1,\ldots,n$.  (Consequently it sends $q_{2i-1}$ to $q_{2i-1}$ for $2i-1\leqslant n$.)
\end{itemize}
All of these facts follow from the constructions of \cite{\Hatcher}
\end{note}

\subsection{Homology of $\Omega_0 SO(n+1)$.}\label{OmegaZeroSubsection}
Let $n\geqslant 2$ and let $\Omega_0 SO(n+1)$ denote the component of $\Omega SO(n+1)$ consisting of contractible loops.  In this subsection we recall Bott's computation of the Hopf algebras $H_\ast(\Omega_0SO(n+1);\Z)$ for $n\geqslant 2$.  We will use these results in the next subsection to describe the Hopf algebras $H_\ast(\Omega SO(n+1);\Q)$ and $H_\ast(\Omega SO(n+1);\Z_2)$.  All references in this section are to \cite[\S\S9-10]{\Bott}.

Bott shows that for $m\geqslant 1$ there are classes
\[\sigma_0,\ldots,\sigma_{2m-1}\in H_\ast(\Omega_0 SO(2m+1);\Q),\quad\sigma_0,\ldots,\sigma_{2m},\varepsilon\in H_\ast(\Omega_0 SO(2m+2);\Q)\]
that generate the Hopf algebras $H_\ast(\Omega_0SO(2m+1);\Q)$, $H_\ast(\Omega_0 SO(2m+2);\Q)$ respectively.   The degrees of these elements are $|\sigma_i|=2i$ and $|\varepsilon|=2m$.  Bott also shows that $\sigma_0-1=0$ and 
\begin{equation}\label{OddRelation}\sigma_i^2-2\sigma_{i-1}\sigma_{i+1}+\cdots\pm 2\sigma_0\sigma_{2i}=0,\qquad i=1,\ldots,m-1\end{equation}
form a complete set of relations among the generators of $H_\ast(\Omega_0 SO(n+1);\Q)$, and that $\sigma_0-1$ and 
\begin{gather}\label{EvenRelationOne}(\sigma_m+\varepsilon)(\sigma_m-\varepsilon)-2\sigma_{m-1}\sigma_{m+1}+\cdots\pm 2\sigma_0\sigma_{2m}=0,\\
\label{EvenRelationTwo}\sigma_i^2-2\sigma_{i-1}\sigma_{i+1}+\cdots\pm 2\sigma_0\sigma_{2i}=0,\qquad i=1,\ldots,m-1
\end{gather}
form a complete set of relations among the generators of $H_\ast(\Omega_0 SO(2m+2);\Q)$.  The coproducts of the generators of $H_\ast(\Omega_0 SO(2m+1);\Q)$ are given by
\begin{equation}
\label{DiagonalZero}D_\ast\sigma_i = {\textstyle\sum}\sigma_{i-j}\otimes\sigma_j\ \ \mathrm{for}\ \ i\leqslant 2m-1
\end{equation}
and the coproducts of the generators of $H_\ast(\Omega_0 SO(2m+2);\Q)$ are given by
\begin{eqnarray}
\label{DiagonalOne}D_\ast\sigma_i &=& {\textstyle\sum}\sigma_{i-j}\otimes\sigma_j\ \ \mathrm{for}\ \ i\leqslant 2m-1,\\
\label{DiagonalTwo}D_\ast\sigma_{2m} &=& {\textstyle\sum}\sigma_{i-j}\otimes\sigma_j + (-1)^m\varepsilon\otimes\varepsilon,\\
\label{DiagonalThree}D_\ast\varepsilon &=& \varepsilon\otimes 1 + 1\otimes\varepsilon.
\end{eqnarray}
This completes the description of the Hopf algebras $H_\ast(\Omega_0 SO(2m+1);\Q)$ and $H_\ast(\Omega_0 SO(2m+2);\Q)$.

Since the homology groups $H_\ast(\Omega_0 SO(n+1);\Z)$ are free over $\Z$ we may regard them as subgroups of $H_\ast(\Omega_0 SO(n+1);\Q)$.  Bott shows that the Hopf algebras $H_\ast(\Omega_0 SO(2m+1);\Z)$ and $H_\ast(\Omega_0 SO(2m+2);\Z)$ are generated by the classes
\[\sigma_0,\ldots,\sigma_{m-1},2\sigma_m,\ldots,2\sigma_{2m-1},\quad
\sigma_0,\ldots,\sigma_{m-1},\sigma_{m}+\varepsilon,\sigma_m-\varepsilon,2\sigma_{m+1},\ldots,2\sigma_{2m}\]
respectively.  The coproducts of these integral generators are completely determined by the coproducts of the rational generators.  The relations among these integral generators are given by $\sigma_0-1=0$ and \eqref{OddRelation} in the case of $H_\ast(\Omega_0SO(2m+1);\Z)$ and $\sigma_0-1=0$ and \eqref{EvenRelationOne}, \eqref{EvenRelationTwo} in the case of $H_\ast(\Omega_0SO(2m+2);\Z)$.

\begin{note}\label{OldLoopHomologyNote}
It follows from Bott's construction that the Hopf algebra morphism $H_\ast(\Omega_0 SO(2m+1);\Q)\to H_\ast(\Omega_0 SO(2m+2);\Q)$ is given by $\sigma_i\mapsto\sigma_i$ for $i=0,\ldots,2m-1$ and that $H_\ast(\Omega_0 SO(2m+2);\Q)\to H_\ast(\Omega_0 SO(2m+3);\Q)$ is given by $\varepsilon\mapsto 0$ and $\sigma_i\mapsto\sigma_i$ for $i=0,\ldots,2m$.
\end{note}

\subsection{Rational and Mod $2$ Homology of $\Omega SO(n+1)$.}\label{OmegaSubsection}

Let $n\geqslant 2$.  We will use the results of the last subsection to describe the Hopf algebras $H_\ast(\Omega SO(n+1);\Q)$ and $H_\ast(\Omega SO(n+1);\Z_2)$.  The main results are summarized as follows:

\begin{proposition}\label{RationalLoopHomologyProposition}
Let $m\geqslant 1$.  Then
\begin{gather*}
H_\ast(\Omega SO(2m+1);\Q)=\Q[\alpha_0,\ldots,\alpha_{2m-1}]/U_m,\\
H_\ast(\Omega SO(2m+2);\Q)=\Q[\alpha_0,\ldots,\alpha_{2m-1},\varepsilon_m]/U_m
\end{gather*}
where the degrees of the generators, and the relations $U_m$, are as described in Theorem~\ref{RationalSOnTheorem}.  The comultiplication is given by
\[D_\ast\alpha_i = \sum\alpha_{i-j}\otimes\alpha_{j},\qquad D_\ast\varepsilon_m = \varepsilon_m\otimes 1 + 1\otimes\varepsilon_m.\]
\end{proposition}

\begin{proposition}\label{ModTwoLoopHomologyProposition}
Let $m\geqslant 1$. Then
\begin{gather*}H_\ast(\Omega SO(2m+1);\Z_2)=\Z_2[a_0,\ldots,a_{m-1},b_0,\ldots,b_{m-1}]\big\slash R_m,\\
H_\ast(\Omega SO(2m+2);\Z_2)=\Z_2[a_0,\ldots,a_{m},b_0,\ldots,b_{m-1}]\big\slash R_m\end{gather*}
where the degrees of the generators, and the relations $R_m$, are as described in Theorem~\ref{ModTwoSOnTheorem}.  The comultiplication is given by
\[D_\ast a_i = \sum a_{i-j}\otimes a_j,\qquad D_\ast b_i = \sum \left(b_{i-j}\otimes a_j + a_j\otimes b_{i-j}\right).\]
\end{proposition}

\begin{proof}[Proof of Proposition~\ref{RationalLoopHomologyProposition}.]
Since $\pi_1SO(n+1)=\Z_2$ for $n\geqslant 2$, we have $H_\ast(\Omega SO(n+1); \Q)\cong \Q[\Z_2]\otimes H_\ast(\Omega_0 SO(n+1); \Q)$ as Hopf algebras.  If we write $\omega$ for the generator of $\Z_2$ then $D_\ast(\omega\otimes 1)=\omega\otimes\omega\otimes 1$.

By the last paragraph and the results of \S\ref{OmegaZeroSubsection}, $H_\ast(\Omega SO(2m+1);\Q)$ is generated by $\omega\otimes 1$ and $1\otimes\sigma_0,\ldots,1\otimes\sigma_{2m-1}$ subject only to the relations arising from $\sigma_0-1=0$, $\omega^2-1=0$, and \eqref{OddRelation}.  The coproducts of these generators are determined by the last paragraph and by \eqref{DiagonalOne}.  The description of $H_\ast(\Omega SO(2m+1);\Q)$ follows by setting $\alpha_i=\omega\otimes\sigma_i$ for $i=0,\ldots,2m-1$.  

Similarly, by the first paragraph and the results of \S\ref{OmegaZeroSubsection}, $H_\ast(\Omega SO(2m+2);\Q)$ is generated by $\omega\otimes 1$ and $1\otimes\sigma_0,\ldots,1\otimes\sigma_{2m},1\otimes\varepsilon$, subject only to the relations arising from $\sigma_0-1=0$, $\omega^2-1=0$, \eqref{EvenRelationOne} and \eqref{EvenRelationTwo}.  The coproducts are determined by the last paragraph and by \eqref{DiagonalTwo} and \eqref{DiagonalThree}.  The description of $H_\ast(\Omega SO(2m+2);\Q)$ follows by setting $\alpha_i=\omega\otimes\sigma_i$ for $i=0,\ldots,2m$ and setting $\varepsilon_m=1\otimes\varepsilon$.
\end{proof}

\begin{proof}[Proof of Proposition~\ref{ModTwoLoopHomologyProposition}.]
As in the last proof, since $\pi_1SO(n+1)=\Z_2$ for $n\geqslant 2$, we have $H_\ast(\Omega SO(n+1); \Z_2)\cong \Z_2[\Z_2]\otimes H_\ast(\Omega_0 SO(n+1); \Z_2)$ as Hopf algebras.  If we write $\omega$ for the generator of $\Z_2$ then $D_\ast(\omega\otimes 1)=\omega\otimes\omega\otimes 1$.

By the last paragraph and the results of \S\ref{OmegaZeroSubsection}, $H_\ast(\Omega SO(2m+1);\Z_2)$ is generated by $\omega\otimes 1$, $1\otimes\sigma_0,\ldots,1\otimes\sigma_{m-1},1\otimes2\sigma_{m},\ldots,1\otimes2\sigma_{2m-1}$, subject to the relations arising from $\sigma_0-1=0$, $\omega^2-1=0$ and \eqref{OddRelation}.  The coproducts of these generators are determined by the last paragraph and by \eqref{DiagonalOne}.   The description of $H_\ast(\Omega SO(2m+1);\Z_2)$ follows by setting $a_i=\omega\otimes\sigma_i$ and $b_i=\omega\otimes2\sigma_{m+i}$ for $i=0,\ldots,m-1$.

Similarly, by the first paragraph and the results of \S\ref{OmegaZeroSubsection}, $H_\ast(\Omega SO(2m+2);\Z_2)$ is generated by $\omega\otimes 1$, $1\otimes\sigma_0,\ldots,1\otimes\sigma_{m-1},1\otimes(\sigma_m\pm\varepsilon), 1\otimes2\sigma_{m+1},\ldots,1\otimes2\sigma_{2m}$, subject to the relations arising from $\sigma_0-1=0$, $\omega^2-1=0$ and \eqref{EvenRelationOne}, \eqref{EvenRelationTwo}.  The coproducts of these generators are determined by the last paragraph and by \eqref{DiagonalTwo} and \eqref{DiagonalThree}.   The description of $H_\ast(\Omega SO(2m+2);\Z_2)$ follows by setting $a_i=\omega\otimes\sigma_i$ and $b_i=\omega\otimes2\sigma_{m+i}$ for $i=0,\ldots,m-1$, and setting $a_m=\omega\otimes(\sigma_m+\varepsilon)$.
\end{proof}

\begin{note}\label{LoopHomologyNote}
Using Note~\ref{OldLoopHomologyNote} we have the following facts.  Let $m\geqslant 1$.
\begin{itemize}
\item The map $H_\ast(SO(2m+1);\Q)\to H_\ast(SO(2m+2);\Q)$ sends $\alpha_i$ to $\alpha_i$ for $i=0,\ldots,2m-1$.
\item The map $H_\ast(SO(2m-1);\Q)\to H_\ast(SO(2m+1);\Q)$ sends $\alpha_i$ to $\alpha_i$ for $i=0,\ldots,2m-3$.
\item The map $H_\ast(SO(2m+1);\Z_2)\to H_\ast(SO(2m+2);\Z_2)$ sends $a_i$ to $a_i$ and $b_i$ to $b_i$ for $i=0,\ldots,m-1$.
\item The map $H_\ast(SO(2m+2);\Z_2)\to H_\ast(SO(2m+3);\Z_2)$ sends $a_i$ to $a_i$ for $i=0,\ldots,m$, sends $b_0$ to $0$, and sends $b_i$ to $b_{i-1}$ for $i=1,\ldots,m-1$.
\item In particular, from the last two facts it follows that for $N$ large enough the maps $H_\ast(SO(2m+1);\Z_2)\to H_\ast(SO(N+1);\Z_2)$, $H_\ast(SO(2m+2);\Z_2)\to H_\ast(SO(N+1);\Z_2)$ send $b_0,\ldots,b_{m-1}$ to $0$.
\end{itemize}
\end{note}

\subsection{The homology suspension.}\label{SuspensionSubsection}
Let $n\geqslant 2$.  The goal of this subsection is to compute the homology suspensions
\begin{gather*}
\sigma\colon H_\ast(\Omega SO(n+1);\Q)\to H_{\ast+1}(SO(n+1);\Q),\\
\sigma\colon H_\ast(\Omega SO(n+1);\Z_2)\to H_{\ast+1}(SO(n+1);\Z_2),
\end{gather*}
using the results of \S\ref{HomologySubsection} and \S\ref{OmegaSubsection}.

\begin{proposition}\label{RationalSuspensionProposition}
There are nonzero rational numbers $\lambda_1,\lambda_2,\ldots$ and $\mu_1,\mu_2,\ldots$ such that:
\begin{enumerate}
\item The homology suspension $\sigma\colon H_\ast(\Omega SO(2m+1);\Q)\to H_{\ast+1}(SO(2m+1);\Q)$ is given by $\sigma(\alpha_{2i-1})=\lambda_i a_{4i-1}$ for $i=1,\ldots,m$ and $\sigma(\alpha_{2i})=0$ for $i=0,\ldots,m-1$.
\item The homology suspension $\sigma\colon H_\ast(\Omega SO(2m+2);\Q)\to H_{\ast+1}(SO(2m+2);\Q)$ is given by $\sigma(\alpha_{2i-1})=\lambda_i a_{4i-1}$ for $i=1,\ldots,m$, $\sigma(\alpha_{2i})=0$ for $i=0,\ldots,m$, and $\sigma(\varepsilon_m)=\mu_mb_{2m+1}$.
\end{enumerate}
\end{proposition}

\begin{note}
It is possible to improve on Proposition~\ref{RationalSuspensionProposition} and show that the constants $\lambda_1,\lambda_2,\ldots$ and $\mu_1,\mu_2,\ldots$ are all equal to 1.
\end{note}

\begin{corollary}\label{RationalDerivationCorollary}
The derivations $\partial_1,\ldots,\partial_m$ of $H_\ast(\Omega SO(2m+1);\Q)$ corresponding to the basis $a_3,\ldots,a_{4m-1}$ of the odd-degree primitive subspace of $H_\ast(SO(2m+1);\Q)$, and the derivations $\partial_1,\ldots,\partial_m,\partial_\varepsilon$ of $H_\ast(\Omega SO(2m+2);\Q)$ corresponding to the basis $a_3,\ldots,a_{4m-1},b_{2m+1}$ of the odd-degree primitive subspace of $H_\ast(SO(2m+2);\Q)$, can be described as follows:
\begin{itemize}
\item $\partial_{i}$ sends $\alpha_{j}$ to $\lambda_i\alpha_{j-2i+1}$ if $j\geqslant 2i-1$ and sends all other generators to $0$.
\item $\partial_\varepsilon$ sends $\varepsilon_m$ to $\mu_m 1$ and sends all other generators to 0.
\end{itemize}
\end{corollary}
\begin{proof}
The derivations, which are defined in Definition~\ref{DerivationsDefinition}, are given explicitly in terms of the coproducts on $H_\ast(\Omega SO(2m+1);\Q)$ and $H_\ast(\Omega SO(2m+2);\Q)$ and the homology suspensions $\sigma\colon H_\ast(\Omega SO(2m+1);\Q)\to H_{\ast+1}(SO(2m+1);\Q)$ and $\sigma\colon H_\ast(\Omega SO(2m+2);\Q)\to H_{\ast+1}(SO(2m+2);\Q)$.  The coproducts are described in Proposition~\ref{RationalLoopHomologyProposition}, and the homology suspensions are described in Proposition~\ref{RationalSuspensionProposition}.  The result follows immediately.
\end{proof}

\begin{proposition}\label{ModTwoSuspensionProposition}
\begin{enumerate}
\item The homology suspension $\sigma\colon H_\ast(\Omega SO(2m+1);\Z_2)\to H_{\ast+1}(SO(2m+1);\Z_2)$ is given by $\sigma(a_i)=q_{2i+1}$ and $\sigma(b_i)=0$ for $i=0,\ldots,m-1$.
\item The homology suspension $\sigma\colon H_\ast(\Omega SO(2m+2);\Z_2)\to H_\ast(SO(2m+2);\Z_2)$ is given by $\sigma(a_i)=q_{2i+1}$ for $i=0,\ldots,m$ and $\sigma(b_i)=0$ for $i=0,\ldots,m-1$.
\end{enumerate}
\end{proposition}

\begin{corollary}\label{ModTwoDerivartionCorollary}
The derivations $\partial_1,\ldots\partial_m$ of $\H_\ast(SO(2m+1);\Z_2)$ corresponding to the basis $q_{1},q_3,\ldots,q_{2m-1}$ of the odd-degree primitive subspace of $H_\ast(SO(2m+1);\Z_2)$, and the derivations $\partial_1,\ldots\partial_{m+1}$ of $\H_\ast(SO(2m+2);\Z_2)$ corresponding to the basis $q_{1},q_3,\ldots,q_{2m+1}$ of the odd-degree primitive subspace of $H_\ast(SO(2m+2);\Z_2)$, are as described in Theorem~\ref{ModTwoSOnTheorem}.
\end{corollary}
\begin{proof}
The derivations, which are defined in Definition~\ref{DerivationsDefinition}, are given explicitly in terms of the coproducts on $H_\ast(\Omega SO(2m+1);\Z_2)$ and $H_\ast(\Omega SO(2m+2);\Z_2)$ and the homology suspensions $\sigma\colon H_\ast(\Omega SO(2m+1);\Z_2)\to H_{\ast+1}(SO(2m+1);\Z_2)$ and $\sigma\colon H_\ast(\Omega SO(2m+2);\Z_2)\to H_{\ast+1}(SO(2m+2);\Z_2)$.  The coproducts are described in Proposition~\ref{ModTwoLoopHomologyProposition} and the homology suspensions are described in Proposition~\ref{ModTwoSuspensionProposition}.  The result follows immediately.
\end{proof}

The rest of this subsection is given to the proofs of Propositions~\ref{RationalSuspensionProposition} and \ref{ModTwoSuspensionProposition}.
We begin with three lemmas that prove special cases of these propositions.  They all follow easily from the Serre spectral sequence of certain fibrations.  Let us recall that for any $n\geqslant 1$ the homology ring $H_\ast(\Omega S^{n+1}; R)$ is $R[u_n]$, where $u_n\in H_n(\Omega S^{n+1};R)$ is a generator that satisfies $\sigma(u_n)=[S^{n+1}]$.

\begin{lemma}\label{RationalSuspensionLemmaOne}
Let $i\geqslant 1$.  Then there is a nonzero $\lambda_i\in\Q$ such that $\sigma(\alpha_{2i-1})=\lambda_i a_{4i-1}$ in $H_{4i-1}(SO(2i+1);\Q)$.
\end{lemma}
\begin{proof}
In this proof all homology groups are taken with rational coefficients.  Let $F_{2i+1}=SO(2i+1)/SO(2i-1)$ denote the space of orthogonal $2$-frames in $\mathbb{R}^{2i+1}$ and let $\pi_i\colon SO(2i+1)\to F_{2i+1}$ denote the projection map.  $F_{2i+1}$ is a rational homology sphere of dimension $4i-1$ with fundamental class $[F_{2i+1}]={\pi_i}_\ast a_{4i-1}$ (see Note~\ref{HomologyNote}), so that $H_\ast(\Omega F_{2i+1})=\Q[v_{4i-2}]$ where $\sigma(v_{4i-2})=[F_{2i+1}]$.

Since the homology groups $H_\ast(\Omega SO(2i-1))$,  $H_\ast(\Omega SO(2i+1))$ and  $H_\ast(\Omega F_{2i+1})$ are all concentrated in even degrees, the Serre spectral sequence of the fibration $\Omega SO(2i-1)\to\Omega SO(2i+2)\to\Omega F_{2i+1}$ collapses at the $E^2$ term and there is a short exact sequence
\[ 0\to H_{4i}(\Omega SO(2i-1))\to H_{4i}(\Omega SO(2i+1))\to H_{4i}(\Omega F_{2i+1})\to 0.\]
By Note~\ref{LoopHomologyNote} the map $H_{4i}(\Omega SO(2i-1))\to H_{4i}(\Omega SO(2i+1))$ has nonzero cokernel generated by $\alpha_{2i-1}$, and so there is some nonzero $\lambda_i\in\Q$ such that ${\pi_i}_\ast(\alpha_{2i-1})=\lambda_i v_{4i-2}$.

By naturality of the homology suspension (Lemma~\ref{SuspensionNaturalityLemma}) and the last two paragraphs we now have ${\pi_i}_\ast\sigma(\alpha_{2i-1})=\lambda_i\pi_\ast a_{4i-1}$, and since $a_{4i-1}$ spans the primitive subspace of $H_{4i-1}(SO(2i+1))$ the result follows.
\end{proof}

\begin{lemma}\label{RationalSuspensionLemmaTwo}
Let $i\geqslant 1$.  Then if $i$ is odd there is a nonzero $\mu_i\in\Q$ such that $\sigma(\varepsilon_i)=\mu_ib_{2i+1}$ in $H_{2i+1}(SO(2i+2);\Q)$, and if $i$ is even there is nonzero $\mu_i\in\Q$ and a $k_i\in\Q$ such that $\sigma(\varepsilon_i)=\mu_i b_{2i+1}+ k_i a_{2i+1}$ in $H_{2i+1}(SO(2i+2);\Q)$.
\end{lemma}
\begin{proof}
In this proof all homology groups are to be taken with rational coefficients.  Let $\rho_i\colon SO(2i+2)\to S^{2i+1}$ denote the projection map.  By Note~\ref{HomologyNote} this map sends $b_{2i+1}$ to $[S^{2i+1}]$ and sends all $a_{4j-1}$ to $0$.  Since the homology groups of $\Omega SO(2i+1)$, $\Omega SO(2i+2)$ and $\Omega S^{2i+1}$ are concentrated in even degrees, the Serre spectral sequence of the fibration $\Omega SO(2i+1)\to\Omega SO(2i+2)\to\Omega S^{2i+1}$ collapses at the $E^2$ term and so we have a short exact sequence
\[0\to H_{2i}(\Omega SO(2i+1))\to H_{2i}(\Omega SO(2i+2))\xrightarrow{{\Omega\rho_i}_\ast} H_{2i}(\Omega S^{2i+1})\to 0.\]
By Note~\ref{LoopHomologyNote} the map $H_{2i+1}(\Omega SO(2i+1))\to H_{2i+1}(\Omega SO(2i+2))$ has nonzero cokernel generated by $\varepsilon_i$, and so there is some nonzero $\mu_i\in\Q$ such that ${\Omega\rho_i}_\ast \varepsilon_i=\mu_iu_{2i}$.  It now follows from naturality of the homology suspension (Lemma~\ref{SuspensionNaturalityLemma}) that ${\rho_i}_\ast\sigma(\varepsilon_i)=\mu_i{\rho_i}_\ast b_{2i+1}$.  If $i$ is even then $b_{2i+1}$ spans the primitive subspace of $H_{2i+1}(SO(2i+2))$, and if $i$ is odd then $b_{2i+1}$ and $a_{2i+1}$ span the primitive subspace of $H_{2i+1}(SO(2i+2))$.  The result follows.
\end{proof}

\begin{lemma}\label{ModTwoSuspensionLemma}
Let $i\geqslant 1$.  Then $\sigma(a_{i})=q_{2i+1}$ in $H_{2i+1}(SO(2i+2);\Z_2)$.
\end{lemma}
\begin{proof}
In this proof all homology groups are taken with coefficients in $\Z_2$.  As in the last proof let $\rho_i\colon SO(2i+2)\to S^{2i+1}$ denote the projection onto the first column.  By Note~\ref{HomologyNote} this sends $q_{2i+1}$ onto $[S^{2i+1}]$.  Since the homology groups of $\Omega SO(2i+1)$, $\Omega SO(2i+2)$ and $\Omega S^{2i+1}$ are concentrated in even degrees, the Serre spectral sequence of the fibration $\Omega SO(2i+1)\to\Omega SO(2i+2)\to\Omega S^{2i+1}$ collapses at the $E^2$ term and so we have a short exact sequence
\[0\to H_{2i}(\Omega SO(2i+1))\to H_{2i}(\Omega SO(2i+2))\xrightarrow{{\Omega\rho_i}_\ast} H_{2i}(\Omega S^{2i+1})\to 0.\]
Recall that $H_\ast(\Omega S^{2i+1})=\Z_2[u_{2i}]$ where $\sigma(u_{2i})=[S^{2i+1}]$.  By Note~\ref{LoopHomologyNote} the map $H_{2i}(\Omega SO(2i+1))\to H_{2i}(\Omega SO(2i+2))$ has nonzero cokernel generated by $a_i$, and so we have ${\Omega\rho_i}_\ast a_i=u_{2i}$.  It now follows from naturality of the homology suspension (Lemma~\ref{SuspensionNaturalityLemma}) that ${\rho_i}_\ast\sigma(a_i)={\rho_i}_\ast q_{2i+1}$.  Since $q_{2i+1}$ spans the primitive subspace of $H_{2i+1}(SO(2i+2))$ the result follows.
\end{proof}

Now that we have proved the three lemmas above, the general results in Proposition~\ref{RationalSuspensionProposition}  and Proposition~\ref{ModTwoSuspensionProposition} will be deduced using the following commutative diagram:
\begin{equation}\label{ProofDiagram}\xymatrix{
H_\ast(\Omega SO(n+1;R) \ar[r]^-{\sigma}\ar[d] & H_{\ast+1}(SO(n+1);R)\ar[d]\\
H_\ast(\Omega SO(N+1);R)\ar[r]_-{\sigma} & H_{\ast+1}(SO(N+1);R)
}\end{equation}
Here the vertical maps are induced from the inclusions $SO({n+1})\to SO({N+1})$ associated to the standard inclusion $\mathbb{R}^{n+1}\to\mathbb{R}^{N+1}$.  Commutativity follows from naturality of the homology suspension (Lemma~\ref{SuspensionNaturalityLemma}).

\begin{proof}[Proof of Proposition~\ref{RationalSuspensionProposition}.]
Throughout the proof we use homology with coefficients in $\Q$.  

We begin with the proof of the first part.  The fact that $\sigma(\alpha_{2i})=0$ for $i=0,\ldots,m-1$ is immediate from the fact that there are no nonzero primitive elements in $H_\ast(SO(2m+1))$ of degree congruent to $1$ modulo $4$.

Let us turn to diagram~\eqref{ProofDiagram} in the case $R=\Q$, $n=2i$, $N=2m$, with $m\geqslant i\geqslant 1$.  The left-hand vertical map sends $\alpha_{2i-1}$ to $\alpha_{2i-1}$ and the right-hand vertical map sends $a_{4i-1}$ to $a_{4i-1}$.  Since $\sigma(\alpha_{2i-1})=\lambda_i a_{4i-1}$ in $H_\ast(SO(2i+1))$ by Lemma~\ref{RationalSuspensionLemmaOne}, it now follows that $\sigma(\alpha_{2i-1})=\lambda_i a_{4i-1}$ in $H_\ast(SO(2m+1))$.  This proves the first part of the proposition.

We now prove the second part of the proposition.  Take diagram~\eqref{ProofDiagram} in the case $R=\Q$, $n=2m$, $N=2m+1$.  The left-hand vertical map sends $\alpha_{2i}$ to $\alpha_{2i}$ and the upper map sends $\alpha_{2i}$ to $0$ by the first part of the proposition.  Thus $\sigma(\alpha_{2i})=0$.  Take diagram~\eqref{ProofDiagram} in the case $R=\Q$, $n=2i$, $N=2m+1$ for $1\leqslant i\leqslant m$.  The left-hand vertical map sends $\alpha_{2i-1}$ to $\alpha_{2i-1}$ and the right-hand vertical map sends $a_{4i-1}$ to $a_{4i-1}$.  Since $\sigma(\alpha_{2i-1})=\lambda_i a_{4i-1}$ in $H_\ast(SO(2i+1))$ by Lemma~\ref{RationalSuspensionLemmaOne}, it now follows that $\sigma(\alpha_{2i-1})=\lambda_i a_{4i-1}$ in $H_\ast(SO(2m+1))$.

It remains to show that $\sigma(\varepsilon_m)=\mu_m b_{2m+1}$.  When $m$ is even this is the statement of Lemma~\ref{RationalSuspensionLemmaTwo}.  When $m$ is odd, Lemma~\ref{RationalSuspensionLemmaTwo} tells us that $\sigma(\varepsilon_m)=\mu_m b_{2m+1}+k_ma_{2m+1}$ for some $k_m\in\Q$.  Taking diagram~\eqref{ProofDiagram} with $n=2m+1$ and $N=2m+2$ and $R=\Q$, we find that the left-hand vertical map sends $\varepsilon_m$ to $0$ by Note~\ref{LoopHomologyNote} and that the right-hand map sends $b_{2m+1}$ to $0$ and $a_{2m+1}$ to $a_{2m+1}$.  It follows that we must have $k_m=0$.  This completes the proof.
\end{proof}

\begin{proof}[Proof of Proposition~\ref{ModTwoSuspensionProposition}.]
Throughout this proof homology groups are to be taken with coefficients in $\Z_2$.  We will prove both parts of the proposition simultaneously.

Note that $a_0$ is represented by a noncontractible loop in $SO(n+1)$, regarded as a point of $\Omega SO(n+1)$, while $q_1$ is represented by the same loop, this time regarded as an element of $H_1(SO(n+1))$.  
The claims $\sigma(a_0)=q_1$ now follow from the definition of the homology suspension.

Now let us turn to diagram \eqref{ProofDiagram} in the case $R=\Z_2$.  The right-hand map in this diagram is always an injection by Note~\ref{HomologyNote}.

Taking $n=2m$ or $2m+1$ and $N$ large enough, the left-hand vertical map sends the classes $b_0,\ldots,b_{m-1}$ to $0$ by Note~\ref{LoopHomologyNote}, so it follows that these classes must vanish under the homology suspension. 

Fix any $1\leqslant i\leqslant m-1$. By taking $n=2i+1$ and $N=2m$, and using the fact that $\sigma(a_i)=q_i$ in $H_{2i+1}(SO(2i+2))$ from Lemma~\ref{ModTwoSuspensionLemma}, we find that $\sigma(a_i)=q_i$ in $H_{2i+1}(SO(2m+1))$.  Similarly, for any $1\leqslant i \leqslant m$ we may take $n=2i+1$, $N=2m+1$ and use the same fact from Lemma~\ref{ModTwoSuspensionLemma} to see that $\sigma(a_i)=q_i$ in $H_{2i+1}(SO(2m+2))$.  This completes the proof.
\end{proof}

\subsection{Proof of Theorems~\ref{RationalSOnTheorem} and \ref{ModTwoSOnTheorem}}\label{ProofSubsection}

\begin{proof}[Proof of Theorem~\ref{RationalSOnTheorem}.]
Proposition~\ref{RationalHomologyProposition} describes the rings $\H_\ast(SO(n+1);\Q)$ and Proposition~\ref{RationalLoopHomologyProposition} describes the rings $H_\ast(\Omega SO(n+1);\Q)$.  The ring-isomorphisms of Theorem~\ref{RationalSOnTheorem} now follow from Theorem~\ref{MainTheorem}.  Proposition~\ref{RationalHomologyProposition} also describes a basis for the odd-degree part of the primitive subspace of $H_\ast(SO(n+1);\Q)$ and describes the effect of the corresponding derivations of $\H_\ast(SO(n+1);\Q)$, while Corollary~\ref{RationalDerivationCorollary} describes the effect of the corresponding derivations of $H_\ast(\Omega SO(n+1);\Q)$ after rescaling.  The description of the BV-operators in Theorem~\ref{RationalSOnTheorem} now follows from Proposition~\ref{DerivationsProposition}, but with the description of the $\partial_i$ and $\partial_\varepsilon$ replaced by the following:
\begin{itemize}
\item $\partial_{i}$ sends $\alpha_{j}$ to $\lambda_i\alpha_{j-2i+1}$ if $j\geqslant 2i-1$ and sends all other generators to $0$.
\item $\partial_\varepsilon$ sends $\varepsilon_m$ to $\mu_m 1$ and sends all other generators to 0.
\end{itemize}
However, by replacing $\partial_i$ with $\partial_i/\lambda_i$, $\delta_i$ with $\lambda_i\delta_i$, and $\beta_{4i-1}$ with $\beta_{4i-1}/\lambda_i$, and by replacing $\partial_\varepsilon$ with $\partial_\varepsilon/\mu_m$, $\delta_\varepsilon$ with $\mu_m\delta_\varepsilon$, and $\gamma_{2m+1}$ with $\gamma_{2m+1}/\mu_m$, we may assume that the description of the BV-operators given in Theorem~\ref{RationalSOnTheorem} holds exactly.  This proves Theorem~\ref{RationalSOnTheorem}.
\end{proof}

\begin{proof}[Proof of Theorem~\ref{ModTwoSOnTheorem}.]
Proposition~\ref{ModTwoHomologyProposition} describes the rings $\H_\ast(SO(n+1);\Z_2)$ and Proposition~\ref{ModTwoLoopHomologyProposition} describes the rings $H_\ast(\Omega SO(n+1);\Z_2)$.  The ring-isomorphisms of Theorem~\ref{ModTwoSOnTheorem} now follow from Theorem~\ref{MainTheorem}.  Proposition~\ref{ModTwoHomologyProposition} also describes a basis for the odd-degree part of the primitive subspace of $H_\ast(SO(n+1);\Z_2)$ and describes the effect of the corresponding derivations of $\H_\ast(SO(n+1);\Z_2)$, while Corollary~\ref{ModTwoLoopHomologyProposition} describes the effect of the corresponding derivations of $H_\ast(\Omega SO(n+1);\Z_2)$.  The description of the BV-operators in Theorem~\ref{ModTwoSOnTheorem} now follows from Proposition~\ref{DerivationsProposition}.  This completes the proof of Theorem~\ref{ModTwoSOnTheorem}.
\end{proof}

\end{document}